\documentclass{amsart}
\usepackage[T1]{fontenc}
\usepackage{ae}
\usepackage{aecompl}
\usepackage{ucs}
\usepackage[utf8x]{inputenc}
\usepackage{color}
\usepackage[enableskew,vcentermath]{youngtab}
\usepackage{amsfonts}
\usepackage{hhline}
\usepackage{amsthm}
\usepackage{amsmath}
\usepackage{amssymb}
\usepackage{hyperref}
\usepackage[all]{pstricks}

\begin{document}
\newtheorem{Le}{Lemma}[section]
\newtheorem{Ko}[Le]{Lemma}
\newtheorem{Sa}[Le]{Theorem}
\newtheorem{pro}[Le]{Proposition}
\newtheorem{Bem}[Le]{Remark}
\newtheorem{Def}[Le]{Definition}
\newtheorem{Bsp}[Le]{Example}
\newtheorem{Con}[Le]{Conjecture}
\numberwithin{equation}{section}

\renewcommand{\l}{\lambda}
\renewcommand{\L}{\lambda}
\newcommand{\bl}{\bar\lambda}
\newcommand{\bn}{\bar\nu}
\newcommand{\A}{\mathcal{A}}
\newcommand{\B}{\mathcal{B}}
\newcommand{\C}{\mathcal{C}}
\newcommand{\D}{\mathcal{D}}
\renewcommand{\a}{\alpha}
\renewcommand{\b}{\beta}
\renewcommand{\c}{\gamma}
\renewcommand{\d}{\delta}
\renewcommand{\k}{\kappa}
\newcommand{\U}{\Upsilon}
\newcommand{\h}{\hfil}
\newcommand{\X}{X}
\renewcommand{\pm}[1]{\begin{pmatrix}#1\end{pmatrix}}
\newcommand{\pinw}{{\pi_{nw}}}
\newcommand{\abs}[1]{\lvert #1 \rvert}
\newcommand{\tm}{\tilde\mu}
\newcommand{\tn}{\tilde\nu}
\newcommand{\lm}{\l/\mu}
\renewcommand{\ln}{\l/\nu}
\newcommand{\ab}{\a/\b}
\newcommand{\hlm}{\hat\l/\mu}
\newcommand{\hab}{\hat\a/\b}
\newcommand{\tlm}{\tilde\l/\mu}
\newcommand{\tab}{\tilde\a/\b}
\newcommand{\m}{\mu}
\newcommand{\n}{\nu}
\renewcommand{\h}{\hfil}

\title[Equality Of Multiplicity Free Skew Characters]{Equality Of Multiplicity Free Skew Characters}
\author[C. Gutschwager]{Christian Gutschwager}
\address{Institut für Algebra, Zahlentheorie und Diskrete Mathematik, Leibniz Universität Hannover,  Welfengarten 1, D-30167 Hannover}
\email{gutschwager (at) math (dot) uni-hannover (dot) de}


\begin{abstract}
In this paper we show that two skew diagrams $\lm$ and $\ab$ can represent the same multiplicity free skew character $[\lm]=[\ab]$ only in the the trivial cases when $\lm$ and $\ab$ are the same up to translation or rotation or if $\l=\a$ is a staircase partition $\l=(l,l-1,\ldots,2,1)$ and $\lm$ and $\ab$ are conjugate of each other.
\end{abstract}
\maketitle

\section{Introduction}
The question under which circumstances two different skew diagrams $\lm$ and $\ab$ give rise to the same skew character $[\lm]=[\ab]$ has lately received much attention and is by \cite[Theorem 4.2]{Gut} equivalent to the question under which circumstances two products of Schubert classes $\sigma_{\a_1}\cdot\sigma_{\a_2}$ and $\sigma_{\b_1}\cdot\sigma_{\b_2}$ are equal.

Trivial cases for equality of skew characters $[\lm]=[\ab]$ are given if the skew diagrams $\lm$ and $\ab$ are the same up to translation or rotation.

In \cite{RSW} and later in \cite{MW} a method for constructing skew diagrams with nontrivial equality of their corresponding skew characters was presented. The fact that for a staircase partition $\l=(l,l-1,\ldots,2,1)$ the skew diagram $\lm$ and its conjugate $(\lm)'$ give rise to the same skew character was proved in \cite[Theorem 7.32]{RSW}.

On the other hand new necessary conditions for two skew diagrams $\lm$ and $\ab$ to give rise to the same skew character have been given recently in \cite{MN} and \cite{HL}.

In \cite{Gut} we classified the skew diagrams $\lm$ whose corresponding skew character $[\lm]$ is multiplicity free which means that in the decomposition $[\lm]=\sum c(\l;\n,\m) [\n]$ all coefficients $c(\l;\n,\m)$ are either $0$ or $1$.

Using algebraic arguments Reiner et al. showed in \cite{RSW} that equality of $[\lm]=[\ab]$ when $\lm$ decays into partitions is only possible in the trivial cases.

We will examine in this paper the case when $\lm$ and $\ab$ are connected skew diagrams and $[\lm]=[\ab]$ is multiplicity free. We will see that the only nontrivial case for $[\lm]=[\ab]$ is when $\l=\a$ is a staircase partition and additionally $\lm=(\ab)'$.

\section{Notation and Littlewood-Richardson-Symmetries}
We mostly follow the standard notation in \cite{Sag} or \cite{Stanley}. A partition $\l=(\l_1,\l_2,\ldots,\l_l)$ is a weakly decreasing sequence of non-negative integers where only finitely many of the $\l_i$ are positive. We regard two partitions as the same if they differ only by the number of trailing zeros and call the positive $\l_i$ the parts of $\l$. The length is the number of positive parts and we write $l(\l)=l$ for the length and $\abs{\l}=\sum_i \l_i$ for the sum of the parts. With a partition $\l$ we associate a diagram, which we also denote by $\l$, containing $\l_i$ left-justified boxes in the $i$-th row and we use matrix-style coordinates to refer to the boxes.

\begin{Def}
We will call a partition $\l$ with $n$ different parts a $dp=n$ partition or write $dp(\l)=n$.
\end{Def}

On the set of partitions we define the lexicographic order in the following way. For two partitions $\l,\m$, we say $\l < \m$ if there is an $i$ with $\l_i < \m_i$ and for $j<i$ we have $\l_j=\m_j$.

The conjugate $\l'$ of $\l$ is the diagram which has $\l_i$ boxes in the $i$-th column.

For $\m \subseteq \l$ we define the skew diagram $\lm$ as the difference of the diagrams $\l$ and $\m$ defined as the difference of the set of the boxes. Rotation of $\lm$ by $180^\circ$ yields a skew diagram $(\lm)^\circ$ which is well defined up to translation. A skew tableau $T$ is a skew diagram in which positive integers are written into the boxes.  We refer with $T(i,j)$ to the entry in box $(i,j)$. A semistandard tableau of shape $\lm$ is a filling of $\lm$ with positive integers such that the following expressions hold for all $(i,j)$ for which they are defined: $T(i,j)<T(i+1,j)$ and $T(i,j)\leq T(i,j+1)$. The content of a semistandard tableau $T$ is $\n=(\n_1,\ldots)$ if the number of occurrences of the entry $i$ in $T$ is $\n_i$. The reverse row word of a tableau $T$ is the sequence obtained by reading the entries of $T$ from right to left and top to bottom starting at the first row. Such a sequence is said to be a lattice word if for all $i,n \geq1$ the number of occurrences of $i$ among the first $n$ terms is at least the number of occurrences of $i+1$ among these terms. The Littlewood-Richardson (LR) coefficient $c(\l;\m,\n)$ equals the number of semistandard tableaux of shape $\lm$ with content $\n$ such that the reverse row word is a lattice word. We will call those tableaux LR tableaux. The LR coefficients play an important role in different contexts (see \cite{Sag} or \cite{Stanley}).

The irreducible characters $[\l]$ of the symmetric group $S_n$ are naturally labeled by partitions $\l\vdash n$. The skew character $[\lm]$ of a skew diagram $\lm$ is defined by the LR coefficients:
\[ [\lm]=\sum_\n c(\l;\m,\n) [\n]. \]

There are many known symmetries of the LR coefficients.

We have that $c(\l;\m,\n)=c(\l;\n,\m)$. The translation symmetry gives $[\lm]=[\ab]$ if the skew diagrams of $\lm$ and $\ab$ are the same up to translation while rotation symmetry gives $[(\lm)^\circ]=[\lm]$. Another well known symmetry is the conjugation symmetry $c(\l';\m',\n')=c(\l;\m,\n)$ .

We say that a skew diagram $\D$ decays into the disconnected skew diagrams $\A$ and $\B$ if no box of $\A$ (viewed as boxes in $\D$) is in the same row or column as a box of $\B$. We write $\D=\A\otimes\B$ if $\D$ decays into $\A$ and $\B$. A skew diagram is connected if it does not decay.

A skew character whose skew diagram $\D$ decays into disconnected (skew) diagrams $\A,\B$ is equivalent to the product of the  characters of the disconnected diagrams induced to a larger symmetric group. We have  \[[\D]=([\A]\times[\B])\uparrow_{S_n\times S_m}^{S_{n+m}}=:[\A]\otimes[\B]\] with $\abs{\A}=n,\abs{\B}=m$.  If $\D=\lm$ and $\A,\B$ are proper partitions $\a,\b$ we have:
\[[\lm]= \sum_\n c(\l;\m,\n)[\n]=\sum_\n c(\n;\a,\b)[\n] =[\a]\otimes[\b].\]

The product of two Schubert classes $\sigma_\a, \sigma_\b$ in the cohomology ring \linebreak $H^*(Gr(l,\mathbb{C}^n),\mathbb{Z})$ of the Grassmannian $Gr(l,\mathbb{C}^n)$ of $l$-di\-men\-sio\-nal subspaces of $\mathbb{C}^n$ is given by:

\[\sigma_\a\cdot\sigma_\b=\sum_{\n\subseteq((n-l)^l)}c(\n;\a,\b)\sigma_\n.\]

In \cite[Section4]{Gut} we established a close connection between the Schubert product and skew characters. To use this relation later on we define the Schubert product for characters in the obvious way as a restriction of the ordinary product:

\[[\a]\star_{(k^l)}[\b]:=\sum_{\n\subseteq(k^l)} c(\n;\a,\b)[\n]. \]

Since translation of $\lm$ does not change the corresponding skew character, we may assume that $\m_i<\l_i$, $\m_i\leq \l_{i+1}$ for each $1\leq i \leq l(\l)$, which means that $\lm$ doesn't have empty rows or columns. We call such a skew diagram a basic skew diagram.

For a basic skew diagram $\lm$ we define two lattice paths from the lower left corner to the upper right corner. The outer lattice path starts to the right, follows the shape of $\l$ and ends upwards in the corner, while the inner lattice path starts upwards, follows the shape of $\m$ and ends with a segment to the right. With $s_{in}$ we refer to the length of the shortest straight segment of the inner lattice path while $s_{out}$ is the length of the shortest straight segment of the outer lattice path.

In \cite[Theorem 3.8]{Gut} we classified the multiplicity free skew characters.
\begin{Sa}[{\cite[Theorem 3.8]{Gut}}]
\label{Sa:mf}
Let $\lm$ be a basic skew diagram which is neither a partition nor a rotated partition. Then $[\lm]$ is multiplicity free if and only if up to rotation of $\lm$ $\m$ is a rectangle and additionally one of the following conditions holds:
\begin{enumerate}
\item $s_{in}=1$
\item $s_{in}=2$ and $\l$ is a $dp=3$ partition
\item $\l$ is a $dp=3$ partition and $s_{out}=1$
\item $\l$ is a $dp=2$ partition.
\end{enumerate}
\end{Sa}

\section{Equality of skew characters}\label{Se:eqcon}
\subsection{Notation and Preliminary}
In this section we prove that two connected skew diagrams $\lm$, $\ab$ which give rise to the same multiplicity free skew character have (up to translation or rotation) to be the same or both $\l$ and $\a$ are the same staircase partition $\l=\a=(n,n-1,n-2,\ldots,2,1)$ and additionally $\m$ is the conjugate of $\b$: $\m=\b'$.

Reiner et al. proved in \cite[Section 6]{RSW} the following, using the Jacobi-Trudi determinant but no LR combinatorics:

\begin{Le}[{\cite[Section 6]{RSW}}]\label{RSW}
Let $\A^1$ and $\A^2$ be skew diagrams and $\A^1$ decay into partitions. Let $[\A^1]=[\A^2]$.

Then the equality is trivial i.e.  up to translation or rotation $\A^1$ and $\A^2$ are the same.

Here translation or rotation does not require the entire skew diagram to be translated or rotated but also includes the case when the partitions into which $\A^1$ decays are translated or rotated independent of each other.
\end{Le}

We will use this lemma always to argue that $\A^1=\c\otimes\d$ with partitions $\c,\d$ and $[\A^1]=[\A^2]$ requires $\A^2=\c\otimes\d$.

To exclude the trivial cases for $[\lm]=[\ab]$ we will in the following always assume that $\m$ is a rectangle and both $\lm$ and $\ab$ are basic skew diagrams.

Furthermore we will use the following notation:

$\l=(\l_1^{l_1},\l_2^{l_2},\ldots)$ with $l=l(\l)=\sum_i l_i$ and $\l_i>\l_{i+1}$, $\m=(\m_1^m)$,

$\a=(\a_1^{a_1},\a_2^{a_2},\ldots)$ with $a=l(\a)=\sum_i a_i$ and $\a_i>\a_{i+1}$, $\b=(\b_1^b)$.

Before we begin with the proofs we state some additional facts.

Recall that the \textit{parts} of a skew diagram $\A=\c/\d$ are the numbers $\c_i - \d_i$ ($1\leq i \leq l(\c)$) and so are the number of boxes in the rows of $\A$. Furthermore the \textit{heights} of a skew diagram $\A$ are the number of boxes in the columns of $\A$ and so are the parts of the conjugated skew diagram. For example the skew diagram  $\A=(7^2,5,3,2)/(4,2^2,1)=${\tiny$\young(::::\hfil\hfil\hfil,::\hfil\hfil\hfil\hfil\hfil,::\hfil\hfil\hfil,:\hfil\hfil,\hfil\hfil)$} has the parts $3,5,3,2,2$ and heights $1,2,3,2,3,2,2$.

It is known, that skew diagrams which give rise to the same skew character have to have the same parts and heights in the same quantity. This follows from the fact, that the LR tableau of the skew diagram $\A$ obtained by filling every column with the entries $1$ to the height of the column ({e.g.\tiny $\young(::::11,:::12,:112,122,233)$}) has as content $\n$ the lexicographic biggest partition whose corresponding character $[\n]$ appears in the decomposition of $[\A]$. Clearly this is the partition obtained by reordering the heights of $\A$ to form a partition. So skew diagrams which give rise to the same skew character have to have the same heights in the same quantity and by conjugation the same holds true for the parts.

In the following we will assume that $[\lm]=[\ab]$. Since we need the same number of rows and columns in $\lm$ and $\ab$ we need $\l_1=\a_1$ and $a=l$ if we assume $[\lm]=[\ab]$.

If we remove from an LR tableau of shape $\lm$ containing $\l_1$ entries $1$ all the boxes with entry $1$ and replace every entry $i>1$ by $i-1$ we obtain an LR tableau of a shape which is obtained from $\lm$ by removing the top box from every column. This gives us a $1-1$ relation between the characters $[\n]\in[\lm]$ with maximal first part $\n_1=\l_1$ and arbitrary characters $[\xi]\in[\hlm]$ with $\hlm$ the skew diagram obtained by removing the top boxes in every column of $\lm$ so $\hat\l=(\l_1^{l_1-1},\l_2^{l_2},\ldots)$ and the $1-1$ relation is given by $\xi_i=\n_{i+1}$. So we have the following lemma:

\begin{Le}\label{Le:remove}
 Let $[\lm]=[\ab]$.

 Let $\hlm$ (resp. $\hab$) be the skew diagram obtained from $\lm$ (resp. $\ab$) by removing the top $i$ boxes from every column of $\lm$ (resp. $\ab$). Let $\tlm$ (resp. $\tab$) be the skew diagram obtained from $\lm$ (resp. $\ab$) by removing in every row of $\lm$ (resp. $\ab$) the left $i$ boxes.

Then $[\hlm]=[\hab]$ and $[\tlm]=[\tab]$.
\end{Le}
\begin{proof}
 $[\hlm]=[\hab]$ follows from the above $1-1$ relation. $[\tlm]=[\tab]$ follows then by conjugation symmetry.
\end{proof}

In most cases we remove the boxes until $\hlm$ decays into two disconnected skew diagrams so that we can use Lemma~\ref{RSW}.

If we say in the following that we remove $l_j$ resp. $\l_i$ boxes from the top resp. left of $\lm$ we always mean that we remove in every column $l_j$ resp. in every row $\l_i$ boxes.

In~\cite{Gut} we proved the following:

\begin{Sa}[{\cite[Theorem 4.2]{Gut}}]
\label{Sa:skewschub}
Let $\m,\l$ be partitions with $\m\subseteq\l\subseteq (k^l)$ with some fixed integers $k,l$. Set $\l^{-1}=(k^l)/\l$.

Then: The coefficient of $[\a]$ in $[\lm]$ equals the coefficient of $[\a^{-1}]=[(k^l)/\a]$ in $[\m]\star_{(k^l)}[\l^{-1}]$.
\end{Sa}

For $k\geq\m_1+\n_1,l\geq l(\m)+l(\n)$ the Schubert product $[\m]\star_{(k^l)}[\n]$ is the ordinary product $[\m]\otimes[\n]$.

\begin{Le}\label{Bem}
Let $[\A^1]=[\A^2]$ with $\A^1$ being an arbitrary skew diagram $\A^1=\c/\d$ having a part $\c_1$ and a height $l(\c)$ (see Figure~\ref{pic:bem}).

Then $\A^1=\A^2$ or $\A^1=(\A^2)^\circ$.
\end{Le}
\begin{proof}
Since $\A^1$ and $\A^2$ have the same heights and parts we have also the part $\c^1$ and the height $l(\c)$ in $\A^2$. The lemma follows now from Lemma~\ref{RSW} and Theorem~\ref{Sa:skewschub}.

\begin{figure}[h]
\psset{xunit=0.35cm,yunit=0.35cm,runit=0.35cm}
\begin{pspicture*}(-4,-0.1)(11,8.1)
\put(-4,4){$\A^1:$}
\pspolygon(0,0)(0,4)(1,4)(1,5)(2,5)(2,6)(3,6)(3,7)(4,7)(4,8)(8,8)(8,3)(7,3)(7,2)(6,2)(6,1)(5,1)(5,0)
\psline{<->}(4.5,0)(4.5,8)\put(5,6){$l(\c)$}
\psline{<->}(0,3.5)(8,3.5)\put(1,2.5){$\c_1$}
\psline[linestyle=dashed](0,0)(0,8)(8,8)(8,0)(0,0)
\psline{->}(-1,7)(2,7)\put(-2,6.5){$\d$}
\psline{->}(9,2)(7,1)\put(9,2){$\c^{-1}$}
\end{pspicture*}
\caption{Lemma~\ref{Bem}: $\A^1$} \label{pic:bem}
\end{figure}

\end{proof}

\subsection{Proofs}
The proof is arranged as follows. We first prove that $[\lm]=[\ab]$ and $dp(\l)=2$ requires $\lm=\ab$ or $\lm=(\ab)^\circ$. Next we prove that $[\lm]=[\ab]$ and $dp(\l)=3$ requires $dp(\a)=3$. Then we prove that $\lm=\ab$ is required if both $\l$ and $\a$ are $dp=3$ partitions and $[\lm]=[\ab]$.

Then we examine the cases when both $\l$ and $\a$ have more than $3$ different parts and $[\lm]=[\ab]$ is multiplicity free.

\begin{Sa}\label{Sa:2parts}
 Let $[\lm]=[\ab]$ and $\l$ be a $dp=2$ partition.

Then $\a$ is also a $dp=2$ partition and $\ab=\lm$ or $\ab=(\lm)^\circ$.
\end{Sa}
\begin{proof}
We will assume that $\a$ has $n\geq 2$ different parts, $dp(\a)=n\geq2$, and show first that $n=2$.

\textit{Case~1: $\l_2>\m_1,\:l_1>m$.} This case is covered by Lemma~\ref{Bem}.

So we have either \textit{Case~2: $l_1\leq m$} or \textit{Case~3: $\l_2\leq\m_1$}. We cannot have both at the same time without $\lm$ decaying into two disconnected rectangles. It is sufficient to examine only the \textit{Case~2} since \textit{Case~3} follows then by conjugation symmetry.

\textit{Case~2: $l_1 \leq m, \l_2 >\m_1$}. By comparing the heights of length $l$ and the parts $\l_1$ in $\lm$ and $\ab$ we get $a_1\leq b, \a_n > \b_1$.

Since $\l_2 > \m_1$ we have $\l_2-\m_1$ times the height $l$ in $\lm$ and so need $\l_2-\m_1=\a_n-\b_1$.

If we remove in $\lm$ the left $\l_2-\m_1$ boxes the remaining skew diagram is: \[\tlm=((\underbrace{\l_1-\m_1-(\l_2-\m_1)}_{=\l_1-\l_2})^{l_1})\otimes((\underbrace{\l_2-\l_2+\m_1}_{=\m_1})^{l-m}).\]

If we remove in $\ab$ the left $\l_2-\m_1=\a_n-\b_1<\a_n$ boxes the remaining skew diagram has to decay by Lemmas~\ref{RSW}  and~\ref{Le:remove} but decays only for $b\geq l-a_n$ and is then:
\[\tab=((\l_1-\a_n)^{a_1},(\a_2-\a_n)^{a_2},\ldots(\a_{n-1}-\a_n)^{a_{n-1}})\otimes((\b_1)^{l-b}) .\]

\begin{figure}[h]
\psset{xunit=0.5cm,yunit=0.5cm,runit=0.5cm}
\begin{pspicture*}(-2,0)(17,5)
\put(-2,2.5){$\tlm:$}
\pspolygon(0,0)(0,2)(3,2)(3,5)(7,5)(7,3)(5,3)(5,0)
\psline(0,0)(2,2)\psline(0,2)(2,0)
\psline(3,5)(5,2)\psline(3,2)(5,5)
\psline[linestyle=dashed](2,0)(2,2)(5,2)(5,5)

\put(8,2.5){$\tab:$}
\pspolygon(10,0)(10,2)(13,2)(13,5)(17,5)(17,4)(16.5,4)(16.5,3.5)(16,3.5)(16,3)(15.5,3)(15.5,2.5)(15,2.5)(15,0)
\psline(10,0)(12,2)\psline(10,2)(12,0)
\psline(13,5)(15,2)\psline(13,2)(15,5)
\psline[linestyle=dashed](12,0)(12,2)(15,2)(15,5)
\end{pspicture*} \caption{\textit{Case~2:} $\tlm$ and $\tab$} \label{fig:2parts.t}
\end{figure}

By Lemmas~\ref{RSW} and~\ref{Le:remove} $\tlm$ and $\tab$ have to be related by translation or rotation if $[\tlm]=[\tab]$ and so we have $[\lm]\neq[\ab]$ for $n\geq 3$.

For $n=2$ we have:
\[\tab=((\l_1-\a_2)^{a_1})\otimes(\b_1^{l-b}).\]

By Lemmas~\ref{RSW} and~\ref{Le:remove} we have either
\begin{align*}
 &\textnormal{\textit{Case 2.1:}}& (\l_1-\l_2)^{l_1}&=(\l_1-\a_2)^{a_1}, & \m_1^{l-m}&=\b_1^{l-b}
\end{align*}
or
\begin{align*}
 &\textnormal{\textit{Case 2.2:}}& (\l_1-\l_2)^{l_1}&=\b_1^{l-b}, & \m_1^{l-m}&=(\l_1-\a_2)^{a_1}.
\end{align*}

In \textit{Case 2.1} we immediately get $\lm=\ab$.

In \textit{Case 2.2} we have
\begin{align*}
\l_1-\l_2&=\b_1,&\l_1-\a_2&=\m_1,&b&=l-l_1,& l-a_1&=m
\end{align*}

which means $\lm=(\ab)^\circ$.
\end{proof}

\begin{Sa}\label{Sa:3parts}
 Let $\l$ be a $dp=3$ partition and $[\lm]=[\ab]$.

Then $dp(\a)=3$.
\end{Sa}
\begin{proof}
By Theorem~\ref{Sa:2parts} $\a$ cannot have only $2$ different parts. So we will assume, that $dp(\a)=n\geq4$ and show that we get contradictions.

\textit{Case 1: $l_1\leq m,\:\l_3\leq\m_1$}. Since we do not have parts $\l_1$ or heights $l$ in $\lm$ we also have $a_1\leq b$ and $\a_n \leq \b_1$.

For $\lm$ to be connected we need $\l_2>\m_1$ and $l_1+l_2>m$.

The cases \textit{Case 1.1: $\l_3\leq \l_2-\m_1$} and \textit{Case 1.2: $l_1\leq l_1+l_2-m$} are related by conjugation symmetry so it is sufficient to consider only \textit{Case 1.1}.

\textit{Case 1.1: $\l_3\leq \l_2-\m_1$}. If we remove in $\lm$ $\l_3$ boxes from the left and in total $l \cdot \l_3$ boxes the remaining skew diagram is $((\l_1-\l_3)^{l_1},(\l_2-\l_3)^{l_2})/\m$ and is not a partition and decays (in the case $\l_2-\m_1=\l_3$) only after $\l_3$ boxes from the left are removed and is connected if only $\l_3-1$ boxes are removed from the left.

If we remove in $\ab$ $\l_3$ boxes from the left the remaining skew diagram $\tab$ has to be by Lemma~\ref{Le:remove} and Theorem~\ref{Sa:2parts} a $dp=2$ partition with a rectangle removed. This cannot be obtained. The best we can obtain is the following. If we demand that $\tab$ is a $dp=2$ partition with or without another partition removed and demand that removing $\l_3$ boxes from the left removes in total $l \cdot \l_3$ boxes we need $n=4,\:\a_4=\l_3$  as well as either $b=a_1, \:\a_1-\b_1-\l_3=0$ or $b=a_1+a_2,\:\a_2-\b_1-\l_3=0$ (see Figure~\ref{fig:3parts.c11}. The case $b=a_1+a_2+a_3,\:\a_3-\b_1-\l_3=0$ can be excluded, because in this case either $\ab$ decays or has a height $l$.). In both cases removing $\l_3-1$ boxes from the left of $\ab$ gives a decaying skew diagram.

\begin{figure}[h]
\psset{xunit=0.4cm,yunit=0.4cm,runit=0.4cm}
\begin{pspicture*}(-1,0)(30,8)
\put(-1,6){$\tlm:$}
\pspolygon(0,0)(0,4)(2,4)(2,8)(7,8)(7,5)(5,5)(5,2.5)(2,2.5)(2,0)
\psline(0,0)(2,4)\psline(0,4)(2,0)
\psline(2,4)(4,8)\psline(2,8)(4,4)
\psline[linestyle=dashed](2,0)(2,4)(4,4)(4,8)

\put(8,4){$\tab^1:$}
\pspolygon(11,0)(11,6)(16,6)(16,8)(18,8)(18,6)(17,6)(17,4)(15,4)(15,2)(13,2)(13,0)
\psline(11,0)(13,6)\psline(11,6)(13,0)
\psline(16,8)(18,6)\psline(16,6)(18,8)
\psline[linestyle=dashed](13,0)(13,6)(18,6)

\put(18,4){$\tab^2:$}
\pspolygon(21,0)(21,4)(24,4)(24,8)(28,8)(28,6)(26,6)(26,4)(25,4)(25,2)(23,2)(23,0)
\psline(21,0)(23,4)\psline(21,4)(23,0)
\psline(24,8)(26,4)\psline(24,4)(26,8)
\psline[linestyle=dashed](23,0)(23,4)(26,4)(26,8)

\end{pspicture*}\caption{\textit{Case 1.1}: $\tlm$ and $\tab^i$}\label{fig:3parts.c11}
\end{figure}

\textit{Case 1.3: $\l_3> \l_2-\m_1,\: l_1>  l_1+l_2-m$}.

If we remove in $\lm$ $\l_2-\m_1$ boxes from the left and in total $l \cdot (\l_2-\m_1)$ boxes the remaining skew diagram $\tlm$ decays into a rectangle and a $dp=2$ partition. If we remove only $\l_2-\m_1-1$ boxes from the left the remaining skew diagram does not decay. We have (see Figure~\ref{fig:3parts.c13lm}):
\begin{equation}
\tlm=((\l_1-\l_2)^{l_1})\otimes (\m_1^{l_1+l_2-m},(\l_3-\l_2+\m_1)^{l_3}). \label{eq:tlm}
\end{equation}

If we remove in $\lm$ $l_1+l_2-m$ boxes from the top which are in total $\l_1 \cdot (l_1+l_2-m)$ boxes the remaining skew diagram $\hlm$ decays also into a rectangle and a $dp=2$ partition. Again removing only $l_1+l_2-m-1$ boxes from the top gives a connected skew diagram. We have (see Figure~\ref{fig:3parts.c13lm}):

\begin{equation}
\hlm=( (\l_1-\m_1)^{m-l_2},(\l_2-\m_1)^{l_2})\otimes (\l_3^{l_3}). \label{eq:hlm}
\end{equation}

\begin{figure}[h]
\psset{xunit=0.5cm,yunit=0.5cm,runit=0.5cm}
\begin{pspicture*}(-2,0)(19,8)
\put(-2,4){$\tlm:$}
\pspolygon(0,0)(0,4)(4,4)(4,8)(8,8)(8,5)(6,5)(6,2)(4,2)(4,0)
\psline(0,0)(2,4)\psline(0,4)(2,0)
\psline(4,4)(6,8)\psline(4,8)(6,4)
\psline[linestyle=dashed](2,0)(2,4)(6,4)(6,8)

\put(9,4){$\hlm:$}
\pspolygon(11,0)(11,4)(15,4)(15,8)(19,8)(19,4)(17,4)(17,2)(14,2)(14,0)
\psline(15,8)(19,6)\psline(19,8)(15,6)
\psline(11,4)(15,2)\psline(15,4)(11,2)
\psline[linestyle=dashed](11,2)(15,2)(15,6)(19,6)

\end{pspicture*}\caption{\textit{Case 1.3}: $\tlm$ and $\hlm$}\label{fig:3parts.c13lm}
\end{figure}

If we remove in $\ab$ $\l_2-\m_1$ boxes from the left the remaining skew diagram has by Lemma~\ref{Le:remove} also to decay into a rectangle and a $dp=2$ partition. To obtain this and remove in total $l \cdot (\l_2-\m_1)$ boxes and also be in the situation that after removing $\l_2-\m_1-1$ boxes from the left the remaining skew diagram $\tab$ does not decay we need $n=4,\: \a_4=\l_2-\m_1$ as well as either $a_1\leq b <a_1+a_2 ,\:\a_2-\b_1= \l_2-\m_1$ or $a_1+a_2\leq b < a_1+a_2+a_3,\:\a_3-\b_1=\l_2-\m_1$ (see Figure~\ref{fig:3parts.c13tab}).
\begin{figure}[h]
\psset{xunit=0.5cm,yunit=0.5cm,runit=0.5cm}
\begin{pspicture*}(-2,0)(19,8)
\put(-2,4){$\tab^1:$}
\pspolygon(0,0)(0,5)(4,5)(4,8)(8,8)(8,6)(6,6)(6,4)(4,4)(4,2)(2,2)(2,0)
\psline(0,0)(2,5)\psline(0,5)(2,0)
\psline(4,5)(6,8)\psline(4,8)(6,5)
\psline[linestyle=dashed](2,0)(2,5)(6,5)(6,8)

\put(8.5,4){$\tab^2:$}
\pspolygon(11,0)(11,3)(13,3)(13,8)(19,8)(19,6)(17,6)(17,4)(15,4)(15,2)(13,2)(13,0)
\psline(11,0)(13,3)\psline(11,3)(13,0)
\psline(13,3)(15,8)\psline(13,8)(15,3)
\psline[linestyle=dashed](13,0)(13,3)(15,3)(15,8)
\end{pspicture*}\caption{\textit{Case 1.3}: The two cases of $\tab^i$}\label{fig:3parts.c13tab}
\end{figure}

With the same reasoning  the skew diagram $\hab$ obtained by removing $l_1+l_2-m$ boxes from the top of $\ab$ must decay into a rectangle and a skew diagram and we need $n=4,\: a_1=l_1+l_2-m$ as well as either $\a_3\leq \b_1 < \a_2, \:a_1+a_2-b=l_1+l_2-m$ or $\a_4\leq \b_1 < \a_3,\: a_1+a_2+a_3-b=l_1+l_2-m$ (see Figure~\ref{fig:3parts.c13hab}).
\begin{figure}[h]
\psset{xunit=0.5cm,yunit=0.5cm,runit=0.5cm}
\begin{pspicture*}(-2,0)(19,8)
\put(-2,4){$\hab^1:$}
\pspolygon(0,0)(0,6)(4,6)(4,8)(8,8)(8,6)(6,6)(6,4)(3,4)(3,2)(2,2)(2,0)
\psline(0,6)(4,4)\psline(0,4)(4,6)
\psline(4,8)(8,6)\psline(4,6)(8,8)
\psline[linestyle=dashed](0,4)(4,4)(4,6)(8,6)

\put(8.5,4){$\hab^2:$}
\pspolygon(11,0)(11,4)(13,4)(13,8)(19,8)(19,6)(17,6)(17,4)(15,4)(15,2)(12,2)(12,0)
\psline(11,2)(13,4)\psline(11,4)(13,2)
\psline(13,8)(19,6)\psline(13,6)(19,8)
\psline[linestyle=dashed](11,2)(13,2)(13,6)(19,6)
\end{pspicture*}\caption{\textit{Case 1.3}: The two cases of $\hab^i$}\label{fig:3parts.c13hab}
\end{figure}

By Lemmas~\ref{RSW} and~\ref{Le:remove} the skew diagrams $\hlm$ resp. $\tlm$ and $\hab$ resp. $\tab$ have to decay into the same partitions. So we get the following:

In the \textit{Case 1.3.1: $a_1\leq b <a_1+a_2,\:\a_2-\b_1= \l_2-\m_1$} we have:
\[\tab^1=((\l_1-\a_2)^{a_1})\otimes((\underbrace{\a_2-\l_2+\m_1}_{=\b_1})^{a_1+a_2-b},(\a_3-\a_2+\b_1)^{a_3}) .\]
Comparing with \eqref{eq:tlm} gives:
\begin{align}
 \a_2&=\l_2, &\b_1&=\m_1 , &\a_3&=\l_3  . \label{eq1}
\end{align}

In the \textit{Case 1.3.2: $a_1+a_2\leq b < a_1+a_2+a_3,\:\a_3-\b_1=\l_2-\m_1$} we have:
\[\tab^2=( (\l_1-\a_3)^{a_1},(\a_2-\l_2-\b_1+\m_1)^{a_2})\otimes(\b_1^{a_1+a_2+a_3-b}).\]
Comparing with \eqref{eq:tlm} gives:
\begin{align}
 \a_3&=\l_1-\m_1 ,   &\a_3-\b_1&=\l_2-\m_1 , & \a_2&=\l_3+\b_1   .\label{eq2}
\end{align}

In the \textit{Case 1.3.x.1: $\a_3\leq \b_1 < \a_2,\: a_1+a_2-b=l_1+l_2-m$} we have:
\[\hab^1= ((\a_2-\b_1)^{a_2}) \otimes (\a_3^{a_3},\a_4^{a_4}) .\]
Comparing with \eqref{eq:hlm} gives:
\begin{align}
 \a_3&\leq \b_1 ,& \a_2&=\l_3+\b_1, & \a_4&=\l_2-\m_1  . \label{eqa}
\end{align}

In the \textit{Case 1.3.x.2: $\a_4\leq \b_1 < \a_3,\: a_1+a_2+a_3-b=l_1+l_2-m$} we have:
\[\hab^2=( (\a_2-\b_1)^{a_2},(\a_3-\b_1)^{a_3})\otimes(\a_4^{a_4}).\]
Comparing with \eqref{eq:hlm} gives:
\begin{align}
 \a_2&=\l_1-\m_1+\b_1  ,&\a_3&=\l_2+\b_1-\m_1,    & \a_4&=\l_3 . \label{eqb}
\end{align}

From the equations of \textit{Case 1.3.1.1} (\eqref{eq1} and \eqref{eqa}) follows $\a_3+\b_1=\l_3+\b_1=\a_2=\l_2=\a_4+\m_1=\a_4+\b_1$, but $\a_3\neq \a_4$.

From the equations of \textit{Case 1.3.1.2} (\eqref{eq1} and \eqref{eqb}) follows $\a_3=\l_3=\a_4$, but $\a_3\neq \a_4$.

From the equations of \textit{Case 1.3.2.1} (\eqref{eq2} and \eqref{eqa}) follows $0\geq \a_3-\b_1=\l_2-\m_1$, but since we are in the case $l_1\leq m$ and $\l_3\leq\m_1$, $\lm$ decays for $\l_2-\m_1\leq 0$.

From the equations of \textit{Case 1.3.2.2} (\eqref{eq2} and \eqref{eqb}) follows $\a_3=\l_1-\m_1=\a_2-\b_1=\l_3=\a_4$, but $\a_4\neq \a_3$.

These contradictions finish \textit{Case 1}.

 \textit{Case 2: $l_1> m,\:\l_3\leq\m_1$} and \textit{Case 3: $l_1\leq m,\:\l_3>\m_1$} are related by conjugation symmetry so it is sufficient to consider only \textit{Case 2}.

\textit{Case 2: $l_1> m,\:\l_3\leq\m_1$}. Since we have the part $\l_1$ $l_1-m$ times in $\lm$ we need $l_1-m=a_1-b$. Removing $l_1-m$ boxes from the top of $\lm$ gives a connected skew diagram $\hlm$ for $\l_2>\m_1$ with
\[\hlm=(\l_1^{m},\l_2^{l_2},\l_3^{l_3})/\m \]
or a decaying skew diagram $\hlm$ for $\l_2\leq \m_1$ with
\[\hlm=((\l_1-\m_1)^{m})\otimes(\l_2^{l_2},\l_3^{l_3}).\]

For $\l_2>\m_1$ removing $l_1-m=a_1-b$ boxes from the top of $\ab$ must yield a connected skew diagram $\hab$ and we get:
\[\hab=(\a_1^{b},\a_2^{a_2},\a_3^{a_3},\a_4^{a_4},\ldots)/\b. \]
Since we are now in \textit{Case 1} with $l_1\leq m,\:\l_3\leq\m_1$ we get from the above $[\lm]\neq[\ab]$.

For $\l_2\leq \m_1$ removing $l_1-m=a_1-b$ boxes from the top of $\ab$ must yield a decaying skew diagram $\hab$ and we get:
\[\hlm=((\a_1-\b_1)^b)\otimes(\a_2^{a_2},\a_3^{a_3},\a_4^{a_4},\ldots).\]
Lemma~\ref{RSW} gives $[\lm]\neq[\ab]$.

\textit{Case 4: $l_1>m,\:\l_3>\m_1$} is covered by Lemma~\ref{Bem}.

\end{proof}

\begin{Sa}\label{Sa:3-3parts}
Let $\l,\a$ be  $dp=3$ partitions and $[\lm]=[\ab]$.

Then $\lm=\ab$.
\end{Sa}\begin{proof}
\textit{Case 1: $l_1\leq m,\: \l_3\leq \m_1$}. Since we do not have heights $l$ and parts $\l_1$ in $\lm$ we also need $a_1\leq b,\: \a_3\leq \b_1$.

For $\lm$ to be connected we need $\l_2>\m_1,\: l_1+l_2>m$.

\textit{Case 1.1: $l_1>l_1+l_2-m$} and \textit{Case 1.2: $\l_3>\l_2-\m_1$} are related by conjugation symmetry so it is sufficient to consider only \textit{Case 1.1}.

\textit{Case 1.1: $l_1>l_1+l_2-m$.}

Removing $l_1+l_2-m$ boxes from the top of $\lm$ gives:

\[\hlm=(\l_3^{l_3})\otimes((\l_1-\m_1)^{m-l_2},(\l_2-\m_1)^{l_2}) .\]

If we remove only $l_1+l_2-m-1$ boxes from the top of $\lm$  the remaining skew diagram does not decay, so we need $l_1+l_2-m=a_1+a_2-b,\: b>a_2$ and removing $a_1+a_2-b$ boxes from the top of $\ab$ gives then:

\[\hab=(\a_3^{a_3})\otimes((\l_1-\b_1)^{b-a_2},(\a_2-\b_1)^{a_2}). \]

This gives $\lm=\ab$.

\textit{Case 1.3: $l_1\leq l_1+l_2-m \Leftrightarrow m\leq l_2 ,\: \l_3\leq \l_2-\m_1$}.

Removing $l_1$ boxes from the top of $\lm$ removes in total $l_1 \cdot \l_1$ boxes and gives a skew diagram which either decays into two disconnected rectangles (for $l_1= l_1+l_2-m$) or is $(\l_2^{l_2},\l_3^{l_3})/\m$ (for $l_1 < l_1+l_2-m$). We need $a_1\leq l_1$ because removing $l_1$ boxes from the top of $\ab$ must either yield a decaying skew diagram or a $dp=2$ partition with a rectangle removed. But since removing the top $l_1$ boxes of $\ab$ has to remove in total $l_1\cdot \l_1$ boxes we need also $a_1\geq l_1$ and so $a_1=l_1$.

Removing $l_1+l_2-m-1$ boxes from the top of $\lm$ yields a connected skew diagram but removing the top $l_1+l_2-m$ boxes yields:
\[\hlm=(\l_3^{l_3})\otimes((\l_2-\m_1)^m) .\]
Since $\ab$ must also decay after removing $l_1+l_2-m$ boxes from the top but not after removing less boxes we get $a_1+a_2-b=l_1+l_2-m$ and so $a_2-b=l_2-m$.

In an analogous way by removing $\l_3$ resp. $\l_2-\m_1$ boxes from the left we get $\a_3=\l_3$ and $\a_2-\b_1=\l_2-\m_1$.

We will first examine the \textit{Case 1.3.1: $\m_1\neq \b_1$} and show that this gives $[\lm]\neq[\ab]$. This covers by conjugation symmetry also the \textit{Case 1.3.2: $m\neq b$}.

For the following construction we only need:
\begin{equation}\label{construction}
\a_2-\b_1=\l_2-\m_1, \qquad a_1=l_1 ,\qquad m\leq l_2, \qquad b \leq a_2.
\end{equation}

Without loss of generality we may also assume that $\m_1 < \b_1$ and set $\b_1=\m_1+n$ which gives $\a_2=\l_2+n$.

\textit{Case 1.3.1.1: $\l_1-(\l_2+n)+1\leq \m_1$}.

We can in $\l$ write entries $1$ to $l_1$ into the columns $\l_2+n$ to $\l_1$ which are in total $\l_1-(\l_2+n)+1$ columns (see Figure~\ref{fig:33parts.1311}). If we place the remaining entries so that they obey the LR rule we get an LR filling of $\l$ with content $\m$ where the box $(l_1,\l_2+n)$ is filled. So there is a character $[\n]$ in $[\lm]$ with $\n$ not containing the box $(l_1,\l_2+n)$.

\begin{figure}[h]
\psset{xunit=0.5cm,yunit=0.5cm,runit=0.5cm}
\begin{pspicture*}(-2,0)(12,8)
\put(-2,4){$\l:$}
\pspolygon(0,0)(0,8)(10,8)(10,5)(6,5)(6,2)(2,2)(2,0)
\put(7.3,7.3){$1$}\put(8.1,7.3){$\cdots$}\put(9.3,7.3){$1$}
\put(7.3,6.3){$\vdots$} \put(9.3,6.3){$\vdots$}
\put(7.3,5.3){$l_1$}\put(8.3,5.3){$\cdots$}\put(9.3,5.3){$l_1$}
\psline(7,5)(7,6)(8,6)(8,5)
\psline{->}(8,4)(7.5,5)
\put(8,3.7){$(l_1,\l_2+n)$}
\end{pspicture*}\caption{LR filling of $\l$ in the \textit{Case 1.3.1.1}}\label{fig:33parts.1311}
\end{figure}

In $\a$ there are $a_2\geq b $ boxes below the box $(l_1,\l_2+n)=(a_1,\a_2)$ so in every LR filling of $\a$ with content $\b$ the box $(l_1,\l_2+n)$ remains empty. So for every character $[\n]$ in $[\ab]$, $\n$ contains a box in position $(l_1,\l_2+n)$. This gives $[\lm]\neq[\ab]$.

\textit{Case 1.3.1.2: $\l_1-(\l_2+n) \geq \m_1 \Leftrightarrow \l_1-\l_2\geq \m_1+n=\b_1$}.

Let us now construct an LR filling of $\a$ with content $\b$ which leaves the box $(l_1+1,\l_2+1)$ empty and so gives $[\lm]\neq[\ab]$.

Place the entries $b$ into the the rows $l$ and $l_1+a_2$ with at least one $b$ in row $l$ (see Figure~\ref{fig:33parts.1312}). Now place for $1< i < b$ the $i$ above $i+1$ if possible. For $b-i=a_3$ we cannot place the entry $i$ above the entry $i+1$ unless there are $\a_2-\a_3$ entries $b$ in row $l_1+a_2$. If there are less than $\a_2-\a_3$ entries $b$ in row $l_1+a_2$ we place the entry $i$ in row $l_1+a_2$ directly left to the $b$. If we now place $\b_1-1$ entries $1$ into the right $\b_1-1$ columns we only get to column $\l_1-(\b_1-1)+1=\l_1-\b_1+2\geq \l_2+2$. If we now place the remaining $1$ atop of one of the entries $2$ which are in row $l_1+a_2+a_3-b+2$ (in one of the columns having in row $a_1+a_2$ an entry $b-a_3$ instead of an entry $b$) then the $1$ is placed in row $l_1+a_2+a_3-b+1\geq l_1+2$ where we used $a_2\geq b$ and $a_3\geq 1$. Also the highest position of a $2$ in this filling is in row $l_1+a_2-b+2$ and so below row $l_1+2$. So the box $(l_1+1,\l_2+1)$ is empty in this LR filling. So we have a character $[\n]$ in $[\ab]$ with $\n$ having a box in position $(l_1+1,\l_2+1)$. But since $(l_1+1,\l_2+1)$ is not in $\l$ there is no character $[\n]$ in $[\lm]$ with $\n$ containing a box in position $(l_1+1,\l_2+1)$ and so $[\lm]\neq[\ab]$.

\begin{figure}[h]
\psset{xunit=0.5cm,yunit=0.5cm,runit=0.5cm}
\begin{pspicture*}(-2,-0.3)(12,8)
\put(-2,4){$\a:$}
\pspolygon(0,-0.3)(0,7.7)(10,7.7)(10,5.7)(8,5.7)(8,1.7)(2,1.7)(2,-0.3)

\put(8,6){$1$}\put(8.5,6){$\cdots$}\put(9.5,6){$1$}

\put(1.3,0){$b$}\put(1.3,0.8){$\vdots$}\put(0.8,0){$b$}\put(0.8,0.8){$\vdots$}
\put(3.5,2){$i$}\put(3.5,2.7){$\vdots$}\put(3.4,3.5){$2$}\put(3.4,4.3){$1$}\put(3,2){$i$}\put(3,2.7){$\vdots$}\put(2.9,3.5){$2$}

\put(4.3,2){$b$}\put(5,2){$\cdots$} \put(6,2){$\cdots$} \put(7.3,2){$b$}
\put(4.3,2.8){$\vdots$}\put(4.3,3.6){$\vdots$}\put(7.3,2.8){$\vdots$}\put(7.3,3.6){$\vdots$}
\put(4.3,4.4){$2$}\put(5,4.4){$\cdots$} \put(6,4.4){$\cdots$} \put(7.3,4.4){$2$}
\put(5.5,5.1){$1$} \put(6,5.1){$\cdots$} \put(7.3,5.1){$1$}

\pspolygon(5.3,5.1)(5.3,5.7)(4.7,5.7)(4.7,5.1)
\psline{->}(6,1)(5,5)

\put(6,0.5){$(l_1+1,\l_2+1)$}

\end{pspicture*}\caption{LR filling of $\a$ in the \textit{Case 1.3.1.2}}\label{fig:33parts.1312}
\end{figure}

\textit{Case 1.3.3: $\m_1=\b_1,\:m=b$}.

Since we have $a_2-b=l_2-m$ and $m=b$ we have also $a_2=l_2$ and from this follows $l_3=l-l_1-l_2=l-a_1-a_2=a_3$.

Also from $\a_2-\b_1=\l_2-\m_1$ and $\b_1=\m_1$ follows $\a_2=\l_2$.

Together with $l_1=a_1$ and $\l_3=\a_3$ which we proved above, this gives the desired $\lm=\ab$ and so finishes the case $l_1\leq m,\: \l_3 \leq \m_1$.

\textit{Case 2: $l_1>m,\: \l_3\leq \m_1$} and \textit{Case 3: $l_1\leq m,\: \l_3>\m_1$} are related by conjugation symmetry so it is sufficient to consider only \textit{Case 2}.

\textit{Case 2: $l_1>m,\: \l_3\leq \m_1$}.

Since we have $l_1-m$ times the part $\l_1$ in $\lm$ we need $l_1-m=a_1-b$.

The skew diagram $\hlm$ obtained after removing $l_1-m$ boxes from the top of $\lm$ decays for $\l_2\leq \m_1$ but is connected for $\l_2 > \m_1$. So we have $\l_2\leq \m_1$ if and only if $\a_2\leq \b_1$.

\textit{Case 2.1: $\l_2\leq \m_1$}.

After removing $l_1-m=a_1-b$ boxes from the top of $\lm$ and $\ab$ we have the skew diagrams

\[\hlm=((\l_1-\m_1)^m) \otimes (\l_2^{l_2},\l_3^{l_3}) \]
and
\[\hab=((\l_1-\b_1)^b)\otimes  (\a_2^{a_2},\a_3^{a_3}) \]

which gives by Lemma~\ref{RSW} $\lm=\ab$.

\textit{Case 2.2: $\l_2> \m_1$}.

Removing $l_1-m=a_1-b$ boxes from the top of $\lm$ and $\ab$ gives
\[\hlm=(\l_1^{m},\l_2^{l_2},\l_3^{l_3})/(\m_1^{m})\]
and
\[\hab=(\l_1^{b},\a_2^{a_2},\a_3^{a_3})/(\b_1^{b}).\]

Using the result of \textit{Case 1:} $l_1\leq m,\: \l_3\leq \m_1$ gives $\ab=\lm$ and finishes \textit{Case 2}.

\textit{Case 4: $l_1>m,\:\l_3>\m_1$} is covered by Lemma~\ref{Bem}.
\end{proof}

We will now compare the multiplicity free skew characters $[\lm]$ and $[\ab]$ when both $\l$ and $\a$ have more than $4$ different parts and assume for the following lemmas that $dp(\l),dp(\a)\geq4$.

\begin{figure}[h]
\psset{xunit=0.35cm,yunit=0.35cm,runit=0.35cm}
\begin{pspicture*}(-2,0)(18.1,10.5)

\pspolygon(0,0)(0,7)(4,7)(4,8)(8,8)(8,7)(7,7)(7,6)(6,6)(6,5)(5,5)(5,4)(4,4)(4,3)(3,3)(3,2)(2,2)(2,1)(1,1)(1,0)
\psline{|-|}(-1,7)(-1,8)
\put(-2,7.1){$1$}
\put(1,5){$m=1$}

\pspolygon(10,0)(10,3)(11,3)(11,8)(18,8)(18,7)(17,7)(17,6)(16,6)(16,5)(15,5)(15,4)(14,4)(14,3)(13,3)(13,2)(12,2)(12,1)(11,1)(11,0)
\psline{|-|}(10,9)(11,9)
\put(10.2,9.3){$1$}
\put(11.5,6){$\m_1=1$}

\end{pspicture*}
\psset{xunit=0.35cm,yunit=0.35cm,runit=0.35cm}
\begin{pspicture*}(-2,0)(20,10.5)
\pspolygon(0,0)(0,6)(5,6)(5,8)(6,8)(6,5)(5,5)(5,4)(4,4)(4,3)(3,3)(3,2)(2,2)(2,1)(1,1)(1,0)
\psline{|-|}(5,9)(6,9)
\put(5.2,9.3){$1$}
\put(-2,7){$\m_1=\l_1-1$}

\pspolygon(7,0)(7,1)(10,1)(10,8)(18,8)(18,7)(17,7)(17,6)(16,6)(16,5)(15,5)(15,4)(14,4)(14,3)(13,3)(13,2)(12,2)(12,1)(11,1)(11,0)
\psline{|-|}(6,0)(6,1)
\put(5,0.1){$1$}
\put(14,1){$m=l-1$}

\end{pspicture*}\caption{The four multiplicity free cases with $\l$ a $dp=n>3$ partition}\label{fig:mf}
\end{figure}
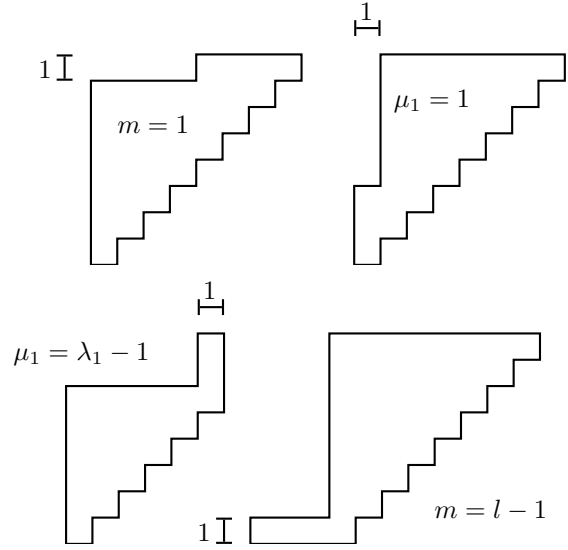

There are $4$ cases when $[\lm]$ is multiplicity free and $\l$ has more than $4$ different parts and we will compare them against each other (see Figure~\ref{fig:mf}).

The first lemma covers by conjugation symmetry also the case when $\m_1=\b_1=1$.
\begin{Le}\label{le:fig}
Let $\lm$ and $\ab$ be skew diagrams with $m=b=1$ and $[\lm]=[\ab]$.

Then $\lm=\ab$.
\end{Le}
\begin{proof}
We have in $\lm$ $l_1-1$ times the part $\l_1$. Since we need in $\ab$ the part $\l_1$ also $a_1-1$ times we get $a_1=l_1$.

If we remove $l_1$ boxes from the top of $\lm$ we either get a connected skew diagram for $\l_2>\m_1,\: l_2>1$ or $\l_3>\m_1$, a disconnected skew diagram for $\l_2>\m_1\geq \l_3,\: l_2=1$ or a partition for $\l_2\leq\m_1$ (see Figure~\ref{fig:mb1}).
Obviously the same must apply for $\hab$ if $[\lm]=[\ab]$.

\begin{figure}[h]
\psset{xunit=0.4cm,yunit=0.4cm,runit=0.4cm}
\begin{pspicture*}(0,0)(28,8)

\pspolygon(0,0)(0,7)(4,7)(4,8)(8,8)(8,6)(6,6)(6,4)(4,4)(4,2)(2,2)(2,0)
\psline(0,7)(4,5)\psline(4,7)(0,5)
\psline(4,8)(8,6)\psline(8,8)(4,6)
\psline[linestyle=dashed](0,5)(4,5)(4,6)(8,6)

\pspolygon(10,0)(10,7)(14,7)(14,8)(18,8)(18,5)(16,5)(16,4)(13,4)(13,2)(12,2)(12,0)
\psline(10,7)(14,4)\psline(14,7)(10,4)
\psline(14,8)(18,5)\psline(18,8)(14,5)
\psline[linestyle=dashed](10,4)(14,4)(14,5)(18,5)

\pspolygon(20,0)(20,7)(25,7)(25,8)(28,8)(28,5)(24,5)(24,4)(23,4)(23,2)(22,2)(22,1)(21,1)(21,0)
\psline(20,7)(24,4)\psline(24,7)(20,4)
\psline(24,7)(25,5)\psline(25,7)(24,5)
\psline(25,8)(28,5)\psline(28,8)(25,5)
\psline[linestyle=dashed](20,4)(24,4)
\end{pspicture*}\caption{Lemma~\ref{le:fig}: The three cases for $\hlm$}\label{fig:mb1}
\end{figure}

Suppose we get a connected skew diagram $\hlm$ with $\hat\l$ having less than $4$ different parts. Then we can use Theorems~\ref{Sa:3parts} and~\ref{Sa:3-3parts} to get $\lm=\ab$.

If $\hlm$ is connected but $\hat\l$ has $4$ or more different parts we can iterate this process until we reach the case where $\hat\l$ has less than $4$ different parts or $\hlm$ is either a disconnected skew diagram or a partition.

Suppose $\hlm$ is a partition.  If we remove only $l_1-1$ boxes from the top of $\lm$ and $\ab$ to get $\hlm^\star$ and $\hab^\star$ we get:
\[\hlm^\star=(\l_1-\m_1)\otimes(\l_2^{l_2},\l_3^{l_3},\l_4^{l_4},\ldots) \]
and
\[\hab^\star=(\l_1-\b_1)\otimes(\a_2^{a_2},\a_3^{a_3},\a_4^{a_4},\ldots) .\]

This gives $\lm=\ab$.

So now suppose $\hlm$ is a disconnected skew diagram. Then we have:

\[\hlm=(\l_2-\m_1)\otimes(\l_3^{l_3},\l_4^{l_4},\ldots)\]
and
\[\hab=(\a_2-\b_1)\otimes(\a_3^{a_3},\a_4^{a_4},\ldots).\]

This gives $\l_i^{l_i}=\a_i^{a_i}$ for $i\geq3$ and $\l_2-\m_1=\a_2-\b_1$.

For $\m_1=\b_1$ we get $\l_2=\a_2$ and since we have in this case also $l_2=a_2=1$ we get $\lm=\ab$.

For $\m_1\neq\b_1$ we can use the construction of LR fillings following equation~\eqref{construction} (page~\pageref{construction})  and get $[\lm]\neq[\ab]$.
\end{proof}

\begin{Le}\label{le:fig2}
Let $\lm$ and $\ab$ be skew diagrams with $\m\neq(1)\neq\b$, $m=1=\b_1$ and $[\lm]=[\ab]$.

Then $\l=(l,l-1,l-2,\ldots,2,1)$ is a staircase partition and $\lm$ is the conjugate of $\ab,\:\ab=\lm'$.
\end{Le}
\begin{proof}
Suppose $l_1>1$. Then we have the part $\l_1$ $l_1-1$ times in $\lm$ and so need $a_1=b+l_1-1$. The smallest height in $\lm$ is either $l_1$ (in the case $\l_2\geq \m_1$) or $l_1-1$ (for $\l_2<\m_1$). The smallest height in $\ab$ is either $a_1=b+l_1-1>l_1$ or $l-b$. For $\l_2\geq \m_1$ this gives $l-b=l_1$  and so $l=b+l_1=a_1+1$ which means that $\a$ has only $2$ different parts. For $\l_2<\m_1$ this gives $l-b=l_1-1$ and so $l=l_1+b-1=a_1$ which means that $\a$ is a rectangle.

So we have $l_1=1$ and since $\lm$ does not decay we need $\l_2>\m_1$.

If we remove the top box of $\lm$ and $\ab$ we get a connected skew diagram for either $l_2>1$ or for $\l_3>\m_1$. Suppose we are in this case then if the new skew diagrams have less than $4$ different parts we can use Theorem~\ref{Sa:3-3parts} which then gives that $[\lm]\neq[\ab]$. If the new skew diagram has $4$ or more different parts we get $\l_2=\a_2$, because $\l_2$ is the number of columns in the new skew diagram, and so $a_1=1$. Since we have the part $\l_1-1$ in $\ab$ we need $\l_2=\l_1-1$, because we need a the part $\l_1-1$ also in $\lm$. We can repeat the above argument until the skew diagram we obtain after removing the top box decays.

So we now assume that the skew diagrams $\hlm$ and $\hab$ obtained after removing the top box of $\lm$ and $\ab$ decay and so we need $l_2=1,\: \l_3\leq\m_1$. If $\hlm$ decays we have:
\[\hlm=(\l_2-\m_1)\otimes(\l_3^{l_3},\l_4^{l_4},\ldots).\]
Since $\hab$ must also decay we need $\a_n=1$ and $a=b+a_n+1$ if $\ab$ has $n$ different parts. We then have:

\[\hab=(1^{a_n})\otimes(\underbrace{(\a_1-1)^{a_1-1}}_{\textnormal{for } a_1>1},(\a_2-1)^{a_2},(\a_3-1)^{a_3},\ldots).\]

Since we have height $1$ in $\lm$ and the smallest height in $\ab$ is either $a_1$ or $l-b=a_n+1>1$ we need $a_1=1$.

Comparing $\hlm$ with $\hab$ gives $\l_2-\m_1=1,\: a_n=1,\: \l_i^{l_i}=(\a_{i-1}-1)^{a_{i-1}}$ for $i=3,\ldots,n$. Since we have again the part $\l_1-1$ in $\ab$ we need $\l_2=\l_1-1$ and so $\m_1=\l_1-2$.

Suppose we have $l_i=a_i=1$ and $\a_i=\l_i$ for $1\leq i < p$ and fixed $p\geq2$. This holds true for $p=2$. Since $l_{i+1}=a_i$ for $i\geq2$ and $l_2=1$ we have also $l_p=1$.

Suppose we have $\l_p>\l_{p+1}+1$ where $\l_{p+1}=0$ is allowed.

We now construct an LR filling of $\l$ with content $\m$. We place in $\l$ entries $1$ into every column of $\l$ but not into the columns $\l_{p+1}+1$ and $\l_{p+1}+2$ and so we get an LR filling which leaves the box $(p,\l_{p+1}+2)$ empty (see Figure~\ref{fig:staircase}). This means that there is a character $[\n]\in[\lm]$ with $\n$ containing a box in position $(p,\l_{p+1}+2)$. Since $\a_p=\l_{p+1}+1$ the box $(p,\l_{p+1}+2)$ is not in $\a$ and so there is no character $[\n]\in[\ab]$ with $\nu$ containing a box in position $(p,\l_{p+1}+2)$. This gives  $[\lm]\neq[\ab]$. So we need $\l_p=\l_{p+1}+1=\a_p$.

\begin{figure}[h]
\psset{xunit=0.4cm,yunit=0.4cm,runit=0.4cm}
\begin{pspicture*}(-2,0)(18,8)

\put(-2,4){$\l:$}
\pspolygon(0,0)(0,8)(8,8)(8,7)(7,7)(7,6)(6,6)(6,5)(5,5)(5,4)(3,4)(3,3)(1,3)(1,0)
\put(0.3,0.3){$1$}\put(1.3,3.3){$1$}\put(2.3,3.3){$1$}\put(5.3,5.3){$1$}\put(6.3,6.3){$1$}\put(7.3,7.3){$1$}
\psline(4,4)(4,5)(5,5)
\psline{->}(3,2)(4.5,4.5) \put(1.5,1){$(p,\l_{p+1}+2)$}

\put(8,4){$\a:$}
\pspolygon(10,0)(10,8)(18,8)(18,7)(17,7)(17,6)(16,6)(16,5)(15,5)(15,3)(14,3)(14,2)(12,2)(12,1)(11,1)(11,0)
\put(10.3,0.3){$6$}\put(11.3,1.3){$5$}\put(13.3,2.3){$4$}\put(15.3,5.3){$3$}\put(16.3,6.3){$2$}\put(17.3,7.3){$1$}
\psline(14,3)(14,4)(15,4)
\psline{<-}(14.5,3.5)(14.5,1.5) \put(13,0.5){$(p+1,\a_p)$}

\end{pspicture*}\caption{Lemma~\ref{le:fig2}: LR fillings of $\l$ and $\a$}\label{fig:staircase}
\end{figure}

So now suppose we have $a_p>1$.

Placing in $\a$ the entries into the rows $1$ to $p-1$ and $p+2$ to $l$ gives an LR filling which leaves the box $(p+1,\a_p)$ empty (see Figure~\ref{fig:staircase}) and so we have a character $[\nu]\in[\ab]$ with $\nu$ containing a box in position $(p+1,\a_p)$. Since $l_p=1$ the $p+1$th row in $\l$ has only $\l_{p+1}<\l_p=\a_p$ boxes which and so there is no character $[\nu]\in[\lm]$ with $\nu$ containing a box in position $(p+1,\a_p)$ and so $[\lm]\neq[\ab]$ for $a_p>1$ and so we need $a_p=1$.

It now follows by induction that $\a=\l$ has to be a staircase partition $\l=(l,l-1,\ldots,2,1)$ which then also gives $b=\m_1$ and so $\lm=(\ab)'$.
\end{proof}

The following lemma covers by conjugation also the case when $\m_1=1,b=a-1=l-1$.
\begin{Le}
Let $\lm$ and $\ab$ be skew diagrams with $m=1,\b_1=\a_1-1=\l_1-1$, $b>1$.

Then $[\lm]\neq[\ab]$.
\end{Le}
\begin{proof}
We have $a_1>b$ since otherwise $\ab$ would decay.

Removing $a_1-b$ boxes from the top of $\ab$ gives $\hab$ which decays into
\[\hab=(1^b)\otimes (\a_2^{a_2},\a_3^{a_3},\a_4^{a_4},\ldots). \]

Since $\m=(\m_1)$ we have that if $\lm$ decays after $a_1-b$ boxes are removed from the top and in total $(a_1-b)\cdot \l_1$ boxes it decays into
\[\hlm=(\l_1-\m_1)\otimes(\l_2^{l_2},\l_3^{l_3},\l_4^{l_4},\ldots). \]

Since $b>1$ we get $[\lm]\neq[\ab]$.
\end{proof}

The following lemma covers by conjugation also the case when $\m_1=1,\b_1=\a_1-1=\l_1-1$.
\begin{Le}
Let $\lm$ and $\ab$ be skew diagrams with $\m_1>1=m, b=l-1=a-1$.

Then $[\lm]\neq[\ab]$.
\end{Le}
\begin{proof}
If we remove the top boxes from $\ab$ we get a partition $\hab$.

If the skew diagram $\hlm$ obtained after removing the top box of every column in $\lm$ is a partition we have $l_1=1$ and $\l_2\leq \m_1$. This means that $\lm$ decays.
\end{proof}

The following lemma covers by conjugation also the case when $m=b=l-1$.
\begin{Le}
Let $\lm$ and $\ab$ be skew diagrams with $\m_1=\b_1=\l_1-1$ and $[\lm]=[\ab]$.

Then $\lm=\ab$.
\end{Le}
\begin{proof}
Since $\lm$ and $\ab$ are connected we need $l_1>m,\:a_1>b$ and because of the part $\l_1$ which therefore exists in $\lm$ and $\ab$ we have $l_1-m=a_1-b$. Removing the top $l_1-m$ boxes of $\lm$ gives:
\[\hlm=(1^m)\otimes(\l_2^{l_2},\l_3^{l_3},\l_4^{l_4}\ldots )\]
and removing the top $l_1-m=a_1-b$ boxes of $\ab$ gives:
\[\hab=(1^b)\otimes(\a_2^{a_2},\a_3^{a_3},\a_4^{a_4}\ldots ).\]

This gives $\lm=\ab$.
\end{proof}

\begin{Le}
Let $\lm$ and $\ab$ be skew diagrams with $m>1,\m_1=\l_1-1, b=l-1$.

Then $[\lm]\neq[\ab]$.
\end{Le}
\begin{proof}
Removing the top box of every column of $\ab$ gives a partition $\hab$.

Since $\m=((\l_1-1)^m)$ with $l_1>m>1$ removing the top boxes of every column of $\lm$ can only give a partition if $\l_2=0$.
\end{proof}

Since the previous $6$ lemmas cover all cases when $[\lm]$ and $[\ab]$ are multiplicity free skew characters with both $\l$ and $\a$ having $4$ or more different parts we get together with the previous theorems the following theorem which holds true without additional prerequisites from this section:

\begin{Sa}
Let $\lm$ and $\ab$ be (connected or decaying) basic skew diagrams and $[\lm]=[\ab]$ multiplicity free.

Then up to translation or rotation
\begin{itemize}
 \item $\lm=\ab$
\item or $\l=\a=(l,l-1,\ldots,2,1)$ and the skew diagrams are conjugate of each other $\lm=(\ab)'$.
\end{itemize}

Here again translation or rotation does not require the entire skew diagram to be translated or rotated, but also allows translation or rotation of the skew diagrams into which $\lm$ decays if it decays.
\end{Sa}

{\bfseries Acknowledgement:} John Stembridge's "SF-package for maple" \cite{stemmaple} was very helpful for computing examples. Furthermore my thanks go to Christine Bessenrodt for helpful discussions and to Stephanie van Willigenburg for pointing out a mistake in a former preprint.

\end{document}